\documentclass[a4paper,twoside,10pt]{article}

\usepackage[a4paper,width=150mm,top=25mm,bottom=25mm,bindingoffset=6mm]{geometry}

\usepackage[english]{babel}

\usepackage{graphicx}
\usepackage[font=small,labelfont=bf]{caption}
\usepackage{subcaption}
\usepackage{color}
\usepackage{mdframed}
\usepackage{framed} 
\usepackage{cancel}
\usepackage[normalem]{ulem}

\usepackage{lineno,hyperref}
\modulolinenumbers[5]

\usepackage{amsmath,amssymb,amsthm,mathtools,mathrsfs,epsfig,amsfonts,wasysym}
\usepackage{MnSymbol}

\usepackage[normalem]{ulem}
\usepackage{ragged2e}
\usepackage{csquotes}
\usepackage{tasks}
\usepackage{paralist}
\usepackage[inline]{enumitem}
\usepackage{tikz-cd}
\usepackage{stmaryrd}

\usepackage[titletoc,title]{appendix}

\newcommand{\e}{\mathrm{e}}
\newcommand{\de}{\mathrm{d}}

\newtheorem{theorem}{Theorem}

\theoremstyle{definition}
\newtheorem{definition}{Definition}

\usepackage{url,hyperref}
\hypersetup{colorlinks=true,linkcolor=blue,filecolor=magenta,urlcolor=cyan}
\usepackage{nicematrix}
\usepackage{tikz}
\usetikzlibrary{fit}

\usepackage[all]{xy}

\newtheorem*{example*}{Example}
\newtheorem*{reminder*}{Reminder}
\newtheorem*{remark*}{Remark}
\newtheorem*{note*}{Note}

\usepackage{authblk}
\usepackage{todonotes}

\definecolor{wildstrawberry}{rgb}{1.0, 0.26, 0.64}

\definecolor{ao(english)}{rgb}{0.0, 0.5, 0.0}

\usepackage{orcidlink}

\begin{document}

\title{A note on the multivariate generalization\\ of a basic simple inequality\vspace{.5cm}\\}

\author[1]{Vasiliki Bitsouni\,\orcidlink{0000-0002-0684-0583}\thanks{\texttt{vbitsouni@math.uoa.gr}}}

\author[1,2]{Nikolaos Gialelis\,\orcidlink{0000-0002-6465-7242}\thanks{\texttt{ngialelis@math.uoa.gr}}}


\affil[1]{Department of Mathematics, National and Kapodistrian University of Athens, Panepistimioupolis, GR-15784 Athens, Greece}

\affil[2]{School of Medicine, National and Kapodistrian University of Athens,\newline GR-11527 Athens, Greece}

\date{}

\maketitle
\begin{center}
\textit{Dedicated to Prof. Ioannis G. Stratis on the occasion of his retirement.\vspace{.5cm}}
\end{center}

\begin{abstract}
\noindent
We introduce the multivariate analogue of the well known inequality $1+x\leq \e^x$, for an abstract non negative real number $x$. The result emerges from the study of the blow up time of certain solutions of the Cauchy problem for a particular ODE. It is also closely related to the notion of completely monotone functions and the theory of divided differences. 
\end{abstract}

\noindent
\textbf{MSC2020-Mathematics Subject Classification System:} 26D07, 34A40

\numberwithin{equation}{section}

\section{Introduction}
\label{intro}

The basic inequality 
\begin{equation}
\label{bsc}
1+x\leq \e^x,
\end{equation}
where $x$ is non negative real number, is common in estimations and useful in applications, especially for (relatively) small values of $x$. 

To the best of the authors' knowledge, there is no multivariate analogue of \eqref{bsc}. The establishment of such a generalization is the aim of this short note. In particular, we show that 
\begin{equation}
\label{gen}
\prod\limits_{i=1}^n{{\left(1+x_i\right)}^{a_i}}\leq\e^{\frac{1}{n}\prod\limits_{i=1}^n{x_i}},
\end{equation}
where $x_1,\dots,x_n$ are pairwise distinct non negative real numbers and $$a_i\coloneqq\frac{\prod\limits_{\substack{j=1\\j\neq i}}^{n}{x_j}}{\prod\limits_{\substack{j=1\\j\neq i}}^{n}{\left(x_j-x_i\right)}}\,.$$ We make the standard empty product convention, hence \eqref{gen} becomes \eqref{bsc} when $n=1$. Moreover, we show that the equality in \eqref{gen} holds only for the case when one $x_i$ equals zero. 

Unlike \eqref{bsc}, \eqref{gen} can not be extended in the whole euclidean space, e.g., when $n=2$, \eqref{gen} is not defined for $\left(x_1,x_2\right)=\left(2,-2\right)$, it is not well defined for $\left(0,-1\right)$, and it does not hold for $\left(1,-\frac{1}{2}\right)$ or $\left(-\frac{1}{4},-\frac{1}{2}\right)$.

The generalized inequality \eqref{gen} naturally emerges from the study of a specific autonomous ODE Cauchy problem. Such an ODE Cauchy problem is fundamental in Mathematical Biology, particularly for the study of one species growth in population dynamics. Knowing \eqref{gen}, a straightforward way can be employed for its demonstration. However, such a way is far from being elementary since it requires both the concept of completely monotone functions and the mean value theorem for divided differences. 

In the present short note, we work as follows: In \hyperref[ode]{Section \ref*{ode}} we employ the ODE approach, i.e. we briefly study the corresponding Cauchy problem in order to show the existence of solutions that blow up in (finite) time (\hyperref[trg]{Section \ref*{trg}}) and we then proceed by proving \eqref{gen}, while we scrutinize the bound of the blow up time of the aforementioned solutions (\hyperref[der]{Section \ref*{der}}). The straightforward path is followed in \hyperref[str]{Section \ref*{str}}, where we first present the preliminaries (\hyperref[prl]{Section \ref*{prl}}), based upon which we obtain the desired inequality (\hyperref[prf]{Section \ref*{prf}}). We conclude our analysis with \hyperref[rpt]{Section \ref*{rpt}}, where we further generalize \eqref{gen} in order to allow repetitions.  

\section{The ODE approach}
\label{ode} 


\subsection{The springboard}
\label{trg} 

We consider the autonomous ODE Cauchy problem   
\begin{equation}
\label{Cpr}
\begin{cases}
\dfrac{\de y}{\de t}{\left(t\right)}=f{\left(y{\left(t\right)}\right)}\coloneqq{\left(-1\right)}^{n+1}y{\left(t\right)}\prod\limits_{i=1}^{n}{\left(1-\frac{1}{k_i}y{\left(t\right)}\right)}\\
y{\left(0\right)}=y_0,\end{cases}
\end{equation}
where $k_1,\dots,k_n$ are pairwise distinct positive real numbers in ascending order. The standard logistic model, and the logistic model with strong Allee effect, see, e.g., \cite{ianpu}, are well known representatives of \eqref{Cpr} for $n=1$ and $n=2$, respectively. 

From the classic theory concerning Cauchy problems for ODEs  (see, e.g., {\cite{hale}}), we can easily deduce existence of a unique smooth maximal solution of \eqref{Cpr}, $$y\colon\,\left(-\varepsilon_1,\varepsilon_2\right)\to\mathbb{R},\text{ where }\varepsilon_{1,2}\in\left(0,\infty\right],$$ that depends smoothly on the problem's data. 

\begin{figure}[!h]
\centering
\includegraphics[width=.75\linewidth]{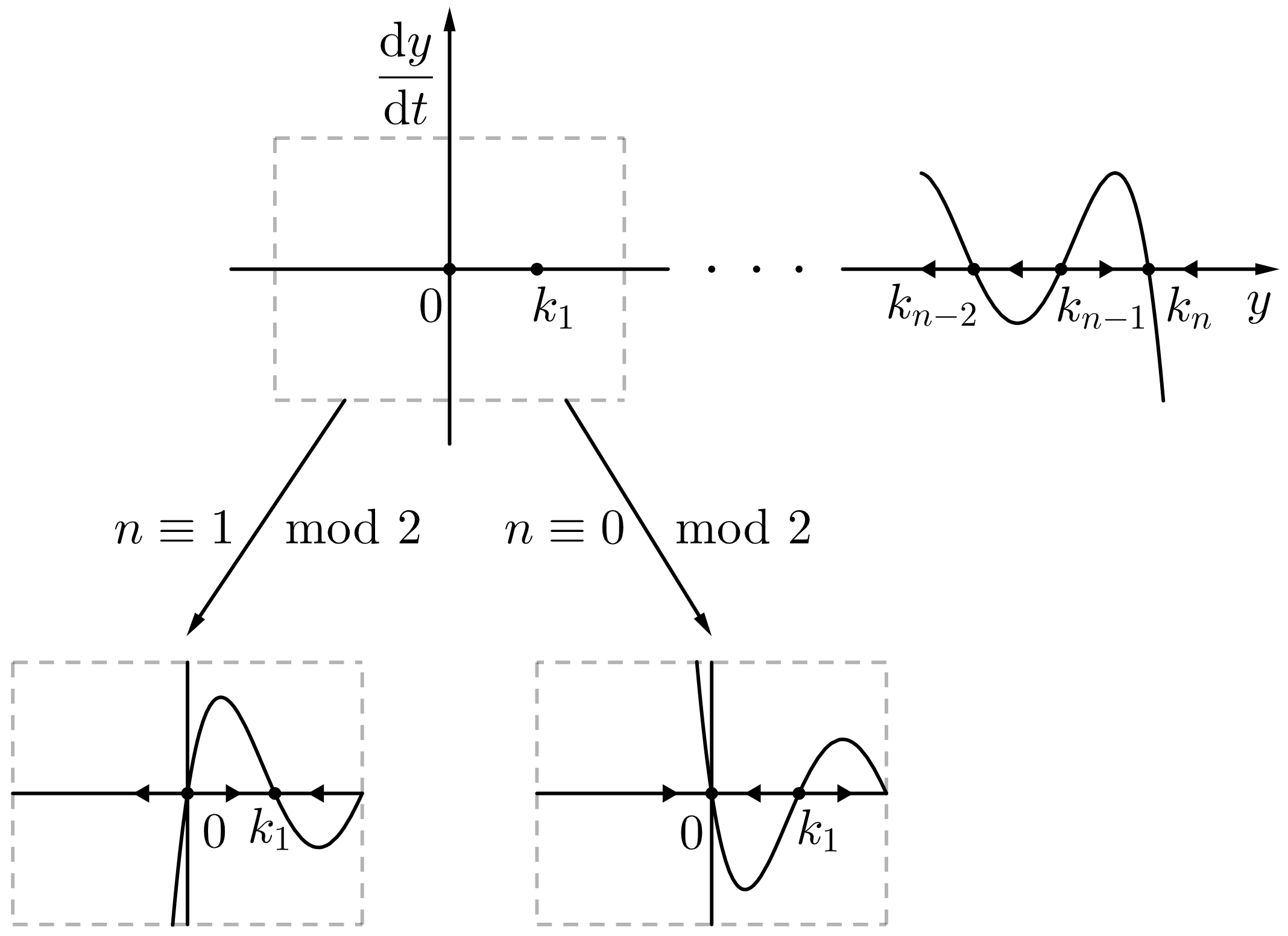}
\caption{Phase line of \eqref{Cpr}.}
\label{phaseline}
\end{figure}

From the phase line, which is depicted in \hyperref[phaseline]{Figure \ref*{phaseline}}, we get that 
\begin{itemize}
\item the sets $$\left(-\infty,0\right),\text{ }\left\{0\right\},\text{ }\left(0,k_1\right),\text{ }\left\{k_1\right\},\text{ }\left(k_1,k_2\right),\dots,\left(k_{n-1},k_n\right),\text{ }\left\{k_n\right\},\text{ }\left(k_n,\infty\right)$$ are (time) invariant, for which we have also utilized the uniqueness of $y$, and 
\item $y$ is $$\begin{cases}
\text{negatively global (in time)},&\text{ when }y_0<0\text{ and }n\equiv 1\mod{2}\\
\text{positively global},&\text{ when }y_0<0\text{ and }n\equiv 0\mod{2}\\
\text{positively global},&\text{ when }y_0>k_n\\
\text{global},&\text{ when }y_0\in\left[0,k_n\right].
\end{cases}$$  
\end{itemize}
It is only left to check the behavior of $y$ when 
\begin{enumerate}
\item $$y_0<0\text{ and }\begin{cases}
t\to{\varepsilon_2}^-,&\text{ when }n\equiv 1\mod{2}\\
t\to{-\varepsilon_1}^+,&\text{ when }n\equiv 0\mod{2}
\end{cases}$$
\item $$y_0>k_n\text{ and }t\to-{\varepsilon_1}^+,$$
\end{enumerate}
in order to close the qualitative analysis of \eqref{Cpr}. 

\begin{enumerate}
\item If $y_0<0$ we split our analysis in two cases: 
\begin{enumerate}
\item If $n\equiv 1\mod{2}$, then $y{\left(t\right)}<0$ and $f{\left(y{\left(t\right)}\right)}<0$ $\forall t\in\left(-\varepsilon_1,\varepsilon_2\right)$, as well as $y$ is monotonically decreasing. We then have 
\begin{multline*}
t=\int\limits_{0}^{t}{dx}=\int\limits_{y{\left(t\right)}}^{y_0}{\left|\frac{1}{f{\left(x\right)}}\right|dx}=\int\limits_{y{\left(t\right)}}^{y_0}{-\frac{1}{f{\left(x\right)}}dx}=\int\limits_{y{\left(t\right)}}^{y_0}{-\frac{1}{{\left(-1\right)}^{n+1}x\prod\limits_{i=1}^{n}{\left(1-\frac{1}{k_i}x\right)}}dx}=\\
=\int\limits_{y{\left(t\right)}}^{y_0}{\frac{\prod\limits_{i=1}^{n}{k_i}}{-x\prod\limits_{i=1}^{n}{\left(k_i-x\right)}}dx},\text{ }\forall t\in\left(0,\varepsilon_2\right).
\end{multline*}
Hence
\begin{equation}
\label{0e1a}
\varepsilon_2=\int\limits_{-\infty}^{y_0}{\frac{\prod\limits_{i=1}^{n}{k_i}}{-x\prod\limits_{i=1}^{n}{\left(k_i-x\right)}}dx}=\int\limits_{-\infty}^{y_0}{\frac{\prod\limits_{i=1}^{n}{k_i}}{{\left(-x\right)}^{n+1}\prod\limits_{i=1}^{n}{\left(1-\frac{k_i}{x}\right)}}dx},
\end{equation}
since $y$ is maximal. Now, from the inequalities 
\begin{equation}
\label{usineq}
\frac{\prod\limits_{i=1}^{n}{k_i}}{{\left(-x\right)}^{n+1}}>0\text{ and }0<\frac{1}{\prod\limits_{i=1}^{n}{\left(1-\frac{k_i}{x}\right)}}<1,\text{ }\forall x<y_0,
\end{equation}
we get the bound 
\begin{equation}
\label{0e1aineq}
\varepsilon_2<\int\limits_{-\infty}^{y_0}{\frac{\prod\limits_{i=1}^{n}{k_i}}{{\left(-x\right)}^{n+1}}dx}=\frac{\prod\limits_{i=1}^{n}{k_i}}{n{\left(-y_0\right)}^{n}}<\infty,
\end{equation}
i.e., $y$ blows up in the future. 
\item If $n\equiv 0\mod{2}$, then $y{\left(t\right)}<0$ and $f{\left(y{\left(t\right)}\right)}>0$ $\forall t\in\left(-\varepsilon_1,\varepsilon_2\right)$, as well as $y$ is monotonically increasing. We then have 
\begin{multline*}
-t=\int\limits_{t}^{0}{dx}=\int\limits_{y{\left(t\right)}}^{y_0}{\left|\frac{1}{f{\left(x\right)}}\right|dx}=\int\limits_{y{\left(t\right)}}^{y_0}{\frac{1}{f{\left(x\right)}}dx}=\int\limits_{y{\left(t\right)}}^{y_0}{\frac{1}{{\left(-1\right)}^{n+1}x\prod\limits_{i=1}^{n}{\left(1-\frac{1}{k_i}x\right)}}dx}=\\
=\int\limits_{y{\left(t\right)}}^{y_0}{\frac{\prod\limits_{i=1}^{n}{k_i}}{-x\prod\limits_{i=1}^{n}{\left(k_i-x\right)}}dx},\text{ }\forall t\in\left(-\varepsilon_1,0\right).
\end{multline*}
Hence
\begin{equation}
\label{0e1b}
\varepsilon_1=\int\limits_{-\infty}^{y_0}{\frac{\prod\limits_{i=1}^{n}{k_i}}{-x\prod\limits_{i=1}^{n}{\left(k_i-x\right)}}dx}=\int\limits_{-\infty}^{y_0}{\frac{\prod\limits_{i=1}^{n}{k_i}}{{\left(-x\right)}^{n+1}\prod\limits_{i=1}^{n}{\left(1-\frac{k_i}{x}\right)}}dx},
\end{equation}
since $y$ is maximal. Inequalities \eqref{usineq} hold here also, thus we get the bound 
\begin{equation}
\label{0e1bineq}
\varepsilon_1<\int\limits_{-\infty}^{y_0}{\frac{\prod\limits_{i=1}^{n}{k_i}}{{\left(-x\right)}^{n+1}}dx}=\frac{\prod\limits_{i=1}^{n}{k_i}}{n{\left(-y_0\right)}^{n}}<\infty,
\end{equation}
i.e. $y$ blows up in the past.
\end{enumerate}
\item If $y_0>k_n$, then $y{\left(t\right)}>k_n$ and $f{\left(y{\left(t\right)}\right)}<0$ $\forall t\in\left(-\varepsilon_1,\varepsilon_2\right)$, as well as $y$ is monotonically decreasing. We then have 
\begin{multline*}
-t=\int\limits_{t}^{0}{dx}=\int\limits_{y_0}^{y{\left(t\right)}}{\left|\frac{1}{f{\left(x\right)}}\right|dx}=\int\limits_{y_0}^{y{\left(t\right)}}{-\frac{1}{f{\left(x\right)}}dx}=\int\limits_{y_0}^{y{\left(t\right)}}{-\frac{1}{{\left(-1\right)}^{n+1}x\prod\limits_{i=1}^{n}{\left(1-\frac{1}{k_i}x\right)}}dx}=\\
=\int\limits_{y_0}^{y{\left(t\right)}}{\frac{\prod\limits_{i=1}^{n}{k_i}}{x\prod\limits_{i=1}^{n}{\left(x-k_i\right)}}dx},\text{ }\forall t\in\left(-\varepsilon_1,0\right).
\end{multline*}
Hence
\begin{equation}
\label{0e2}
\varepsilon_1=\int\limits_{y_0}^{\infty}{\frac{\prod\limits_{i=1}^{n}{k_i}}{x\prod\limits_{i=1}^{n}{\left(x-k_i\right)}}dx}=\int\limits_{y_0}^{\infty}{\frac{\prod\limits_{i=1}^{n}{k_i}}{{x}^{n+1}\prod\limits_{i=1}^{n}{\left(1-\frac{k_i}{x}\right)}}dx},
\end{equation}
since $y$ is maximal. Now, from the inequalities 
\begin{equation*}
\frac{\prod\limits_{i=1}^{n}{k_i}}{{x}^{n+1}}>0\text{ and }0<\frac{1}{\prod\limits_{i=1}^{n}{\left(1-\frac{k_i}{x}\right)}}<\frac{1}{\prod\limits_{i=1}^{n}{\left(1-\frac{k_i}{y_0}\right)}},\text{ }\forall x>y_0,
\end{equation*}
we get the bound 
\begin{equation}
\label{0e2ineq}
\varepsilon_1<\int\limits_{y_0}^{\infty}{\frac{\prod\limits_{i=1}^{n}{k_i}}{{x}^{n+1}\prod\limits_{i=1}^{n}{\left(1-\frac{k_i}{y_0}\right)}}dx}=\frac{\prod\limits_{i=1}^{n}{k_i}}{n{y_0}^{n}\prod\limits_{i=1}^{n}{\left(1-\frac{k_i}{y_0}\right)}}<\infty,
\end{equation}
i.e. $y$ blows up in the past.
\end{enumerate}


\subsection{Derivation of the main result}
\label{der}

With \eqref{0e1a}, \eqref{0e1b} and \eqref{0e2} at hand, we proceed to the calculation of the exact blow up time of the pair of families of extreme solutions of \eqref{Cpr}, i.e. for $y_0<0$ and $y_0>k_n$. 

\textit{First}, we notice that the functions inside the integrals in the aforementioned relations need to be further analyzed: 
\begin{enumerate}
\item We employ the partial fraction decomposition to write 
\begin{equation}
\label{pfm1}
\frac{\prod\limits_{i=1}^{n}{k_i}}{-x\prod\limits_{i=1}^{n}{\left(k_i-x\right)}}=\frac{A}{-x}+\sum\limits_{i=1}^{n}{\frac{A_i}{k_i-x}},\text{ }\forall x<y_0<0.
\end{equation}
For the calculation of the coefficients $A,A_1,\dots,A_n$ we study the consequent equality 
\begin{equation}
\label{cnsq1}
\prod\limits_{i=1}^{n}{k_i}=A\prod\limits_{i=1}^{n}{\left(k_i-x\right)}-x\sum\limits_{i=1}^{n}{A_i\prod\limits_{\substack{j=1\\j\neq i}}^{n}{\left(k_j-x\right)}},\text{ }\forall x\in\mathbb{R}.
\end{equation}
We set $x=0$ in \eqref{cnsq1} to get $A=1$. Moreover, we eliminate the coefficient of $x^n$ in the right hand side of \eqref{cnsq1} to have $${\left(-1\right)}^n\left(A+\sum\limits_{i=1}^{n}{A_i}\right)=0\Rightarrow\sum\limits_{i=1}^{n}{A_i}=-1.$$ Then, we fix an $i_0\in\left\{1,\dots,n\right\}$ and we set $x=k_{i_0}$ in \eqref{cnsq1} to get $$\prod\limits_{i=1}^{n}{k_i}=-k_{i_0}\sum\limits_{i=1}^{n}{A_i\prod\limits_{\substack{j=1\\j\neq i}}^{n}{\left(k_j-k_{i_0}\right)}}=-k_{i_0}A_{i_0}\prod\limits_{\substack{j=1\\j\neq i_0}}^{n}{\left(k_j-k_{i_0}\right)}\Rightarrow A_i=-\frac{\prod\limits_{\substack{j=1\\ j\neq i}}^{n}{k_j}}{\prod\limits_{\substack{j=1\\ j\neq i}}^{n}{\left(k_j-k_i\right)}}.$$ 
\item In an analogous manner, we write 
\begin{equation}
\label{pfm2}
\frac{\prod\limits_{i=1}^{n}{k_i}}{x\prod\limits_{i=1}^{n}{\left(x-k_i\right)}}=\frac{B}{x}+\sum\limits_{i=1}^{n}{\frac{B_i}{x-k_i}},\text{ }\forall x>y_0>k_n,
\end{equation}
and its consequent equality $$\prod\limits_{i=1}^{n}{k_i}=B\prod\limits_{i=1}^{n}{\left(x-k_i\right)}+x\sum\limits_{i=1}^{n}{B_i\prod\limits_{\substack{j=1\\j\neq i}}^{n}{\left(x-k_j\right)}},\text{ }\forall x\in\mathbb{R}.$$ Dealing as above, we get $$B={\left(-1\right)}^n,\text{ }B_i=\frac{\prod\limits_{\substack{j=1\\ j\neq i}}^{n}{k_j}}{\prod\limits_{\substack{j=1\\ j\neq i}}^{n}{\left(k_i-k_j\right)}}={\left(-1\right)}^nA_i\text{ and }\sum\limits_{i=1}^n{B_i}={\left(-1\right)}^n\sum\limits_{i=1}^n{A_i}={\left(-1\right)}^{n+1}.$$ 
\end{enumerate}

\textit{Second}, we employ \eqref{pfm1} and \eqref{pfm2} to calculate $\varepsilon_2$ or $\varepsilon_1$ and $\varepsilon_1$, respectively: 
\begin{enumerate}
\item If $y_0<0$, we then have
\begin{multline*}
\int\limits_{-\infty}^{y_0}{\frac{\prod\limits_{i=1}^{n}{k_i}}{-x\prod\limits_{i=1}^{n}{\left(k_i-x\right)}}dx}=\int\limits_{-\infty}^{y_0}{\frac{1}{-x}+\sum\limits_{i=1}^{n}{\frac{A_i}{k_i-x}}dx}=\int\limits_{-y_0}^{\infty}{\frac{1}{x}+\sum\limits_{i=1}^{n}{\frac{A_i}{k_i+x}}dx}=\\
=\lim\limits_{x\to\infty}{\ln{\left(x\prod\limits_{i=1}^n{{\left(k_i+x\right)}^{A_i}}\right)}}-\ln{\left(-y_0\prod\limits_{i=1}^n{{\left(k_i-y_0\right)}^{A_i}}\right)}.
\end{multline*}
Besides, $$\ln{\left(x\prod\limits_{i=1}^n{{\left(k_i+x\right)}^{A_i}}\right)}=\ln{\left(x^{1+\sum\limits_{i=1}^n{A_i}}\prod\limits_{i=1}^n{{\left(\frac{k_i}{x}+1\right)}^{A_i}}\right)}=\ln{\left(\prod\limits_{i=1}^n{{\left(\frac{k_i}{x}+1\right)}^{A_i}}\right)},$$ hence $$\int\limits_{-\infty}^{y_0}{\frac{\prod\limits_{i=1}^{n}{k_i}}{-x\prod\limits_{i=1}^{n}{\left(k_i-x\right)}}dx}=\ln{\left(\prod\limits_{i=1}^n{{\left(1+\frac{k_i}{-y_0}\right)}^{-A_i}}\right)}$$ and so 
\begin{equation}
\label{vare1}
\ln{\left(\prod\limits_{i=1}^n{{\left(1+\frac{k_i}{-y_0}\right)}^{-A_i}}\right)}=\begin{cases}
\varepsilon_2,&\text{ if }n\equiv 1\mod{2}\\
\varepsilon_1,&\text{ if }n\equiv 0\mod{2}.
\end{cases}
\end{equation}
\item Analogously, if $y_0>k_n$, then
\begin{multline*}
\int\limits_{y_0}^{\infty}{\frac{\prod\limits_{i=1}^{n}{k_i}}{x\prod\limits_{i=1}^{n}{\left(x-k_i\right)}}dx}=\int\limits_{y_0}^{\infty}{\frac{{\left(-1\right)}^n}{x}+\sum\limits_{i=1}^{n}{\frac{B_i}{x-k_i}}dx}=\\
=\lim\limits_{x\to\infty}{\ln{\left(x^{{\left(-1\right)}^n}\prod\limits_{i=1}^n{{\left(x-k_i\right)}^{B_i}}\right)}}-\ln{\left({y_0}^{{\left(-1\right)}^n}\prod\limits_{i=1}^n{{\left(y_0-k_i\right)}^{B_i}}\right)}. 
\end{multline*}
Since $$\ln{\left(x^{{\left(-1\right)}^n}\prod\limits_{i=1}^n{{\left(x-k_i\right)}^{B_i}}\right)}=\ln{\left(x^{{\left(-1\right)}^n+\sum\limits_{i=1}^n{B_i}}\prod\limits_{i=1}^n{{\left(1-\frac{k_i}{x}\right)}^{B_i}}\right)}=\ln{\left(\prod\limits_{i=1}^n{{\left(1-\frac{k_i}{x}\right)}^{B_i}}\right)},$$ we have 
\begin{equation}
\label{vare2}
\varepsilon_1=\ln{\left(\prod\limits_{i=1}^n{{\left(1-\frac{k_i}{y_0}\right)}^{-B_i}}\right)}.
\end{equation}
\end{enumerate}

Now, all we have to do is to combine \eqref{0e1aineq} and \eqref{0e1bineq} with \eqref{vare1}, as well as \eqref{0e2ineq} with \eqref{vare2}, to derive \eqref{gen}:
\begin{enumerate}
\item For $y_0<0$ we have $$\ln{\left(\prod\limits_{i=1}^n{{\left(1+\frac{k_i}{-y_0}\right)}^{-A_i}}\right)}<\frac{\prod\limits_{i=1}^{n}{k_i}}{n{\left(-y_0\right)}^{n}}$$ and we set $x_i\coloneqq\frac{k_i}{-y_0}$ to get the desired result when $x_1,\dots,x_n$ are positive.
\item For $y_0>k_n$ we have $$\ln{\left(\prod\limits_{i=1}^n{{\left(1-\frac{k_i}{y_0}\right)}^{-B_i}}\right)}<\frac{\prod\limits_{i=1}^{n}{k_i}}{n{y_0}^{n}\prod\limits_{i=1}^{n}{\left(1-\frac{k_i}{y_0}\right)}}$$ and now we set $$x_i\coloneqq\frac{\frac{k_i}{y_0}}{1-\frac{k_i}{y_0}}$$ to get the desired result also when $x_1,\dots,x_n$ are positive. 
\item We notice that the equality in \eqref{gen} holds when one $x_i$ equals zero. Since the strict inequality holds for positive values of all $x_i$, we deduce that the equality in \eqref{gen} holds \textit{only} when one $x_i$ equals zero. 
\end{enumerate}

\section{The straightforward approach}
\label{str} 


\subsection{Preliminaries}
\label{prl} 

First we give the definition of completely monotone functions (see, e.g., \cite{miller2001completely}). 
\begin{definition}
\label{cmf}
A function $f\in C^\infty{\left(\left(0,\infty\right)\right)}$ is (strictly) completely monotone iff $${\left(-1\right)}^n f^{\left(n\right)}{\left(x\right)}\geq 0\text{ }\left(\,{\left(-1\right)}^n f^{\left(n\right)}{\left(x\right)}> 0\,\right),\text{ }\forall \left(x,n\right)\,\in\,\left(0,\infty\right)\times\mathbb{N}.$$ 
\end{definition} 
A known example of strictly completely monotone function is 
\begin{equation}
\label{ln}
f{\left(x\right)}\coloneqq \frac{\ln{\left(1+x\right)}}{x},\text{ }\forall x\in\left(0,\infty\right),
\end{equation}
since (see, e.g., \cite{miller2001completely}) $${\left(-1\right)}^{n}f^{\left(n\right)}{\left(x\right)}=n!\int\limits_0^1{\frac{t^n}{{\left(1+tx\right)}^{n+1}}dt}>0,\text{ }\forall \left(x,n\right)\,\in\,\left(0,\infty\right)\times\mathbb{N}.$$ We note that from the classic Beppo Levi theorem of monotone convergence, or the Lebesgue (or Arz\'{e}la) theorem of dominated convergence, we get 
\begin{equation}
\label{nln}
{\left(-1\right)}^n\lim\limits_{x\to 0^+}f^{\left(n\right)}{\left(x\right)}=n!\int\limits_0^1{t^ndt}=\frac{n!}{n+1}.
\end{equation}

We also need a generalization of a well known result to higher derivatives, the mean value theorem for divided differences (see, e.g., \cite{popoviciu1933quelques}, \cite{sahoo1998mean}, or \cite{abel2004mean}). 
\begin{theorem}
\label{cauchypop}
Let $x_1,\dots,x_n$ be pairwise distinct real numbers, with $$m\coloneqq \min\limits_{i\in\left\{1,\dots,n\right\}}{\left\{x_i\right\}}\text{ and }M\coloneqq\max\limits_{i\in\left\{1,\dots,n\right\}}{\left\{x_i\right\}},$$ as well as $f\in \,C{\left(\left[m,M\right]\right)\cap C^{n-1}{\left(\left(m,M\right)\right)}}$. Then $\exists x_0\in\left(m,M\right)$, such that $$\left[x_1,\dots,x_n;f\right]\coloneqq\sum\limits_{i=1}^n{\frac{f{\left(x_i\right)}}{\prod\limits_{\substack{j=1\\j\neq i}}^{n}{\left(x_i-x_j\right)}}}=\frac{f^{\left(n-1\right)}{\left(x_0\right)}}{\left(n-1\right)!}.$$
\end{theorem}

\subsection{Proof of the main result}
\label{prf} 

We assume that $x_1,\dots,x_n$ are pairwise distinct positive real numbers and we want to deduce that $$\prod\limits_{i=1}^n{{\left(1+x_i\right)}^{a_i}}<\e^{\frac{1}{n}\prod\limits_{i=1}^n{x_i}},$$ or, equivalently, 
\begin{equation}
\label{ggenn} 
\sum\limits_{i=1}^n{\frac{f{\left(x_i\right)}}{\prod\limits_{\substack{j=1\\j\neq i}}^{n}{\left(x_j-x_i\right)}}}\,<\,\frac{1}{n}\,,
\end{equation}
where $f$ is as in \eqref{ln}. Employing \hyperref[cauchypop]{Theorem \ref*{cauchypop}}, we have that $\exists x_0\in\left(m,M\right)$, such that $$\sum\limits_{i=1}^n{\frac{f{\left(x_i\right)}}{\prod\limits_{\substack{j=1\\j\neq i}}^{n}{\left(x_j-x_i\right)}}}={\left(-1\right)}^{n-1}\left[x_1,\dots,x_n;f\right]=\frac{{\left(-1\right)}^{n-1}f^{\left(n-1\right)}{\left(x_0\right)}}{\left(n-1\right)!}.$$ Since $f$ is strictly completely monotone, the function ${\left(-1\right)}^nf^{\left(n-1\right)}$ is strictly increasing, thus ${\left(-1\right)}^{n-1}f^{\left(n-1\right)}$ is strictly decreasing, which implies that $$\frac{{\left(-1\right)}^{n-1}f^{\left(n-1\right)}{\left(x_0\right)}}{\left(n-1\right)!}<\frac{{\left(-1\right)}^{n-1}}{{\left(n-1\right)!}}\lim\limits_{x\to 0^+}f^{\left(n-1\right)}{\left(x\right)}\,\overset{\text{\eqref{nln}}}{=}\,\frac{1}{n}$$ and \eqref{ggenn} then follows. 

\section{Allowing repetitions}
\label{rpt}

In \eqref{gen} the numbers $x_1,\dots,x_n$ are distinct. Here we show how to deal with probable repetitions. An elegant approach relies on a proper scaling of \eqref{gen}. Indeed, for some natural numbers $r_1,\dots,r_n$ at hand, we directly deduce from \eqref{gen} that $$\prod\limits_{i=1}^n{\prod\limits_{j=1}^{r_i}{{\left(1+\tau_{ij}x_i\right)}}^{a_{ij}}}\leq\e^{\frac{1}{m}\prod\limits_{i=1}^n{{x_i}^{r_i}\prod\limits_{j=1}^{r_i}{\tau_{ij}}}},$$ where $x_1,\dots,x_n,\tau_{11},\dots,\tau_{nr_n}$ are real numbers, such that $\tau_{11}x_1,\tau_{12}x_1,\dots,\tau_{n\,r_n-1}x_{n},\tau_{nr_n}x_n$ are pairwise distinct non negative real numbers, as well as $$m\coloneqq\sum\limits_{i=1}^{n}{r_i}\text{ and }a_{ij}\coloneqq\frac{\left(\prod\limits_{\substack{k=1\\k\neq i}}^{n}{{\tau_k}^{r_k}\prod\limits_{\ell=1}^{r_k}{\tau_{k\ell}}}\right)\left(\prod\limits_{\substack{\ell=1\\ \ell\neq j}}^{r_i}{\tau_{i\ell}}\right)}{\left(\prod\limits_{\substack{k=1\\k\neq i}}^{n}{\prod\limits_{\ell=1}^{r_k}{\left(\tau_{k\ell}x_k-\tau_{ij}x_i\right)}}\right)\left(\prod\limits_{\substack{\ell=1\\ \ell\neq j}}^{r_i}{\left(\tau_{i\ell}-\tau_{ij}\right)}\right)}\,.$$ Now, we make the least complex choice of numbers, namely we choose $x_1,\dots,x_n$ to be pairwise distinct non negative real numbers and $\tau_{ij}=j$, and the above inequality then becomes 
\begin{equation}
\label{gennn}
\prod\limits_{i=1}^n{\prod\limits_{j=1}^{r_i}{{\left(1+jx_i\right)}}^{a_{ij}}}\leq\e^{\frac{1}{m}\prod\limits_{i=1}^n{r_i!{x_i}^{r_i}}},
\end{equation}
where $x_1,\dots,x_n$ are pairwise distinct non negative real numbers, as well as $$m\coloneqq\sum\limits_{i=1}^{n}{r_i}\text{ and }a_{ij}\coloneqq\frac{\prod\limits_{\substack{\ell=1\\ \ell\neq j}}^{r_i}{\ell}\,\prod\limits_{\substack{k=1\\k\neq i}}^{n}{r_k!{x_k}^{r_k}}}{\prod\limits_{\substack{\ell=1\\ \ell\neq j}}^{r_i}{\left(\ell-j\right)}\,\prod\limits_{\substack{k=1\\k\neq i}}^{n}{\prod\limits_{\ell=1}^{r_k}{\left(\ell x_k-jx_i\right)}}}\,.$$ The equality \eqref{gennn} holds only when one $x_i$ equals zero. 

\section*{Acknowledgement}

The authors are grateful to Dr. Dan \c{S}tefan Marinescu for bringing, both the notion of completely monotone functions (as well as that of the Bernstein functions) and the theory of divided differences, to their attention. 

\end{document}